\input amstex
\documentstyle{amsppt}
\NoRunningHeads
\pagewidth{14cm}
\pageheight{21.1cm}
\magnification=1200
\loadmsbm

\define\End{\operatorname {End}}
\define\red{{\operatorname {red}}}

\define\tr{\operatorname {tr}}

\define\rank{\operatorname {rank}}

\define\spa{\operatorname {span}}

\define\ol{\overline}
\define\vvp{\varphi}

\NoBlackBoxes


\author Yoji Yoshii
\endauthor

\affil Akita National College of Technology \\  1-1 Iijima Bunkyocho
Akita-shi,  Japan 011-8511\\
yoshii\@akita-nct.jp
\endaffil

\topmatter
\title Locally extended affine root systems
\endtitle

\abstract 
We consider a natural generalization of both locally finite irreducible root systems
and extended affine root systems defined by Saito.
We classify the systems.
\endabstract

\footnotetext""{2000 Mathematics Subject Classification 17B20, 17B65, 17B67, 20K15.}
\endtopmatter

\document

\head
Introduction
\endhead

Let us recall the definition of a finite irreducible root system
in a euclidean space $V$,
i.e., $V\thickapprox\Bbb R^n$
with a positive definite form $(\cdot,\cdot)$.

\definition{Definition 0} 
A subset $\frak R$ of $V$ 
is called a {\it finite irreducible root system}
if
\roster
\item"(A1)"  
$0\notin\frak R$ and $\frak R$ spans $V$;

\item"(A2)"
$\langle\alpha,\beta\rangle\in\Bbb Z$ for all 
$\alpha, \beta\in\frak R$,
where $\langle\alpha,\beta\rangle=\frac{2(\alpha,\beta)}{(\beta,\beta)}$;

\item"(A3)"
$\sigma_\alpha(\beta)\in \frak R$ for all $\alpha,\beta\in \frak R$, where
$\sigma_\alpha(\beta)=\beta-\langle\beta,\alpha\rangle\alpha$;

\item"(A4)" 
$\frak R=\frak R_1\cup \frak R_2$ and $(\frak R_1,\frak R_2)=0$ imply $\frak R_1=\emptyset$ or
$\frak R_2=\emptyset$.
($\frak R$ is irreducible.)
\endroster

\enddefinition

We note that $\frak R$ becomes automatically a finite set (see [LN1, 4.2]
or [MY1, Prop. 4.2]).
Needless to say,
these interesting subsets were crucial in the classification of finite-dimensional simple Lie algebras
and of finite reflection groups in the 20th century.
In 1985, K. Saito introduced the notion of a generalized root system [S].
He changed the frame $V$ from the euclidean space to a general vector space
over $\Bbb R$
equipped with a symmetric bilinear form, not necessarily a positive definite form,
and replaced the axiom (A1) to:
$$\text{$(\alpha,\alpha)\neq 0$  for all $\alpha\in\frak R$, and
$\frak R$ spans $V$.}
$$
This change is natural since
$(\alpha,\alpha)\neq 0$ whenever $\alpha\neq 0$
in a euclidean space.
Moreover,
Saito added two extra axioms:
\roster
\item"(A5)"  
the additive subgroup generated by $\frak R$ is a full lattice in $V$;

\item"(A6)"
the codimension of the radical of $V$ is finite.
\endroster
He called such a root system an {\it extended affine root system}
if the form is positive semidefinite.
(Later the notion of an extended affine root system was used in a different sense [A-P],
but it was proved in [A2] there is a natural correspondence between both notions.)
If the dimension of the radical of the positive semidefinite form is 1,
the extended affine root systems
are irreducible affine root systems
in the sense of Macdonald [M].
One of the Saito's main purposes was to construct a Lie algebra
whose anisotropic roots form an extended affine root system
having a 2-dimensional radical.

Our interest now is not Saito's root systems but extended affine root systems.
We generalize Saito's axioms of extended affine root systems with good reasons.
First of all, we make our new concept contain the so-called {\it locally finite irreducible root systems}
(see [LN1]),
which are obtained simply by changing the frame $V$ in Definition 0
to an infinite-dimensional euclidean space,
i.e., an infinite-dimensional vector space over $\Bbb R$
with a positive definite form.  
(Then $\frak R$ becomes automatically a locally finite set,
i.e., $|W\cap\frak R|<\infty$ for any finite-dimensional subspace $W$ of $V$.)
It turns out that if we simply take off the axiom (A6),
locally finite irreducible root systems are contained.
Next, we assume the base field to be $\Bbb Q$, not $\Bbb R$.
Notice that in the setting of
 finite-dimensional simple Lie algebras,
root systems naturally appear in vector spaces over $\Bbb Q$,
and then one gets a euclidean space by simply tensoring with $\Bbb R$.
Besides, our theory of extended affine Lie algebras
also produces a root system in a vector space over $\Bbb Q$ (see [MY1]).
Once we start with a vector space over $\Bbb Q$,
the axiom (A5) is equivalent to saying that
the abelian group  
generated by $\frak R$, say $\langle \frak R\rangle$, is free.
Thus, it seems better to have as axiom that
\roster
\item"(A7)"
$\langle \frak R\rangle$ is free
\endroster
in our setup.
However, we can say much about the classification
without assuming (A7).
So we simply take off the axiom (A5)
(and we do not assume (A7) either),
and  we get our definition of a {\it locally extended affine root system}
in Definition 1.
As a special case,
we call a locally extended affine root system
an {\it extended affine root system} 
if  (A6) and (A7) hold.
Thus our extended affine root systems are the same as  Saito's
if we consider the embedding of ours into the real
vector space $\Bbb R\otimes_\Bbb Q V$.

We classify locally extended affine root systems 
in terms of triples of  reflection spaces
by the methods from [A-P]
(see Theorem 7),
which was also done in [LN2] in a more general setting.
Also, we show some relations between 
the isomorphisms of locally extended affine root systems and
the similarities of
 reflection spaces in Theorem 10.
Then, when $\dim_\Bbb Q V^0=1$, we get more information
by a simple observation about subgroups of $\Bbb Q$ in Corollary 13.
Finally, we give some interesting examples of Lie algebras whose root systems
are locally extended affine root systems.

We thank Professors Saeid Azam, Jun Morita and Erhard Neher for their several suggestions.

\head
Basic Concepts
\endhead

\definition{Definition 1} 
Let $V$ be a vector space over $\Bbb Q$
with a positive semidefinite bilinear form $(\cdot,\cdot)$. 
A subset $\frak R$ of $V$ 
is called a {\it locally extended affine root system}
or a {\it LEARS} for short if
\roster
\item"(A1)" $(\alpha,\alpha)\neq 0$  for all $\alpha\in\frak R$, and
$\frak R$ spans $V$;

\item"(A2)"
$\langle\alpha,\beta\rangle\in\Bbb Z$ for all 
$\alpha, \beta\in\frak R$,
where $\langle\alpha,\beta\rangle=\frac{2(\alpha,\beta)}{(\beta,\beta)}$;

\item"(A3)"
$\sigma_\alpha(\beta)\in \frak R$ for all $\alpha,\beta\in \frak R$, where
$\sigma_\alpha(\beta)=\beta-\langle\beta,\alpha\rangle\alpha$;

\item"(A4)" 
$\frak R=\frak R_1\cup \frak R_2$ and $(\frak R_1,\frak R_2)=0$ imply $\frak R_1=\emptyset$ or
$\frak R_2=\emptyset$.
($\frak R$ is irreducible.)
\endroster

A LEARS $\frak R$ is called {\it reduced} if
$2\alpha\notin \frak R$ for all $\alpha\in\frak R$.

\enddefinition

Note that if $V$ is finite-dimensional and $(\cdot,\cdot)$ is positive definite,
then $\frak R$ is exactly a finite irreducible root system
(see [MY1, Prop. 4.2]).

\medskip

Let $$V^0:=\{x\in V\mid (x,y)=0\ \text{for all $y\in V$}\}$$ be the radical of the form.
Note that
$$V^0=\{x\in V\mid (x,x)=0\}.$$
We call $\dim_\Bbb Q V^0$ the {\it null dimension} of $\frak R$,
which can be any cardinality.

We denote the additive subgroup of $V$ generated by a subset $S$ of $V$
by $\langle S\rangle$.

We call a LEARS
$(\frak R,V)$ an {\it extended affine root system}
or an {\it EARS} for short,
if 
$$\text{$\dim_\Bbb QV/V^0<\infty$ and
$\langle \frak R\rangle$ is free.}
$$
This coincides with the concept, which was firstly introduced by Saito in 1985 [S].
As we mentioned in the Introduction,
the notion of an EARS was also used in a different sense in [A-P],
but Azam showed that there is a natural correspondence between the two notions in [A2].
We use here the Saito's one since he is the first person who defined it and his root system naturally 
generalized Macdonald's {\it affine root systems} in [M].

In Corollary 5 later,
we will see that
if the abelian group
 $\langle \frak R\rangle\cap V^0$ is free, then
$\langle \frak R\rangle$ 
is free.
So the condition that $\langle \frak R\rangle$ is free
can be replaced by the condition that
 $\langle \frak R\rangle\cap V^0$ is free.
 
 Recall the notion of {\it rank} for a torsion-free abelian group $G$,
 that is,
 $\rank G=\dim_\Bbb Q (\Bbb Q\otimes_\Bbb Z G)$.
It is easy to check that 
if $G$
is a subgroup of a $\Bbb Q$-vector space $W$,
then
$\rank G=\dim_\Bbb Q\spa_\Bbb Q G$,
where $\spa_\Bbb Q G$ is the subspace of $W$ spanned by $G$ over $\Bbb Q$.

Thus, in our root system $(\frak R, V)$,
we have
$$\rank (\langle \frak R\rangle\cap V^0)
=\dim_\Bbb Q\spa_\Bbb Q(\langle \frak R\rangle\cap V^0)
=\dim_\Bbb QV^0
=\text{(the null dimension of $\frak R$)}.$$

Now,
when our torsion-free abelian group $\langle \frak R\rangle\cap V^0$ happens to be free,
we say that $\frak R$ has {\it nullity}. 
(We simply want to distinguish the easier case ``free''.)
For example,
$\frak R$ has nullity 1 means that 
$\langle \frak R\rangle\cap V^0\cong \Bbb Z$,
and
$\frak R$ has null dimension 1 means that 
$\langle \frak R\rangle\cap V^0$ is 
isomorphic to a nonzero subgroup of $\Bbb Q$.
Also, by Corollary 5,
if an EARS $\frak R$ has finite nullity,
then $\langle \frak R\rangle$ 
is free of finite rank.
Thus we simply say that $\frak R$ is an EARS of finite rank when the EARS has finite nullity.

Our LEARS are a natural generalization of the existing concept EARS.
In fact, Saito's EARS are the same as
our EARS
embedded into the real vector space $\Bbb R\otimes_\Bbb Q V$.
Similarly, irreducible affine root systems in the sense of Macdonald [M]
are our EARS of nullity 1.
Note that the  
reduced irreducible affine root systems are the real roots of affine Kac-Moody Lie algebras.
The elliptic root systems defined by Saito [S]
are our EARS of nullity 2.
Also, the sets of nonisotropic roots of 
EARS in [A-P] are our reduced EARS of finite rank
(see [A2]).

\medskip

Finally, we call a LEARS of nullity 1
a {\it locally affine root system} or a {\it LARS} for short.

\medskip

Let $(\frak R,V)$ be a LEARS, and $(\bar \frak R,\bar V)$ the canonical image onto $V/V^0$. 
Then $\bar V$ admits the induced positive definite form, and  thus
$$
\text{$(\bar \frak R,\bar V)$ is a locally finite irreducible root system.}
$$
Note that our definition of a locally finite irreducible root system
is a LEARS (in Definition 1) so that the form is positive definite,
and then one can show that the system is in fact locally finite (see 
[LN1, 4.2] or [MY1, Prop.4.2]).

\head
Reflectable bases
\endhead

Locally finite irreducible root systems which are not finite were classified as the reduced types
$\text A_{\frak I}$, $\text B_{\frak I}$, $\text C_{\frak I}$, $\text D_{\frak I}$,
and the nonreduced type $\text{BC}_{\frak I}$ for any infinite index set $\frak I$
(see [LN1, Ch.8]).
More precisely, let $\{\epsilon_i\}_{i\in\frak I}$ be an orthonormal basis for
an infinite-dimensional euclidean space $V$
(or 
an infinite-dimensional vector space over $\Bbb Q$
with positive definite form),
and let
$$
\align
\text A_{\frak I}&=\{\epsilon_i-\epsilon_j\mid i\neq j\in \frak I\},\\
\text B_{\frak I}&=\{\pm\epsilon_i,\ \pm(\epsilon_i\pm\epsilon_j)\mid i\neq j\in \frak I\},\\
\text C_{\frak I}&=\{\pm(\epsilon_i\pm\epsilon_j),\ \pm2\epsilon_i\mid i\neq j\in \frak I\},\\
\text D_{\frak I}&=\{\pm(\epsilon_i\pm\epsilon_j)\mid i\neq j\in \frak I\},\\
\text{BC}_{\frak I}&
=\text B_{\frak I}\cup\text C_{\frak I}.
\endalign
$$
Note that each root system spans $V$ except $\text A_{\frak I}$.
If $|\frak I|=\ell$ is finite, then 
an ordinary notation of the root system is $\text A_{\ell-1}$ instead of $\text A_{\ell}$.
So it might be better to write something like $\text A_{\frak I-1}$
or $\dot{\text A}_{\frak I}$
instead of just $\text A_{\frak I}$.
However, to simplify the notation, we stipulate to write $\text A_{\frak I}$ when
$|\frak I|$ is infinite.

Let $(\frak R,V)$ be a locally finite irreducible root system (including the finite case)
and assume it is reduced.
A basis $\Pi$ of $V$ is called a {\it reflectable base} of $\frak R$
if $\Pi\subset\frak R$ and for any $\alpha\in\frak R$,
$$\alpha=\sigma_{\alpha_1}\cdots\sigma_{\alpha_k}(\alpha_{k+1})$$
for some $\alpha_1, \ldots, \alpha_{k+1}\in\Pi$.
(Any root can be obtained by reflecting a root of $\Pi$ relative to hyperplanes determined by $\Pi$.)
This is a well-known property which a root base of a reduced finite root system possesses.
It is known that a locally finite irreducible root system which is countable possesses a root base,
but this is not the case for a locally finite irreducible root system which is uncountable.
(See [LN1, \S 6].  They also prove that
 there always exists an integral base even in the uncountable case.
 However, it is easy to see that an integral base is not necessarily a reflectable base.)
Thus we need to show the existence of a reflectable base 
in a reduced locally finite irreducible root system which is uncountable, and we have:

\proclaim{Lemma 2}
There exists a reflectable base in any reduced locally finite irreducible root system [LN2, Lem. 5.1].
In particular, the additive subgroup generated by each locally finite irreducible root system is free
(see also [LN1, Thm 7.5]).
\endproclaim

\head
Classification
\endhead

We devote this section to classifying LEARS.
(The argument below is a special case of [LN2, 4.9, 5.2].)

Let $( \frak R,V)$ be a LEARS.
Let $V'$ be a subspace of $V$ so that $V=V'\oplus V^0$,
and
$$\Delta=\Delta_{V'}:=\{\alpha\in V'\mid \alpha+ s\in \frak R\ \text{for some $ s\in V^0$}\}.$$
We note that $\Delta$ is bijectively mapped onto $\bar\frak R$ by $\bar{\phantom a}$.
Moreover, $\bar{\phantom a}$ is a linear isomorphism from $V'$ onto $\bar V$
satisfying $(v', w')=(\bar v', \bar w')$ for all $v',w'\in V'$.
Hence, $(\Delta,V')$ is a locally finite irreducible root system isomorphic to
$(\bar\frak R,\bar V)$.
We often say that $\frak R$ has type $\Delta$.
For each $\alpha\in\Delta$, we set
$$S_{\alpha}:=\{ s\in V^0\mid\alpha+ s\in\frak R\}.$$ 
Then
$$\frak R=\bigsqcup_{\alpha\in\Delta}\ (\alpha+S_{\alpha}).$$ 
Since $\frak R$ spans $V$,
$$
\text{$\cup_{\alpha\in\Delta}\ S_{\alpha}$ spans $V^0$}.
\tag S0
$$
Also,
for any $\alpha+s, \beta+s'$ with $\alpha,\beta\in\Delta$, $s\in S_\alpha$ and $s'\in S_\beta$,
we have
$$
\align
\sigma_{\alpha+s}(\beta+s')
&=\beta+s'-\langle\beta+s',\alpha+s\rangle(\alpha+s)\\
&=\sigma_{\alpha}(\beta)+s'-\langle\beta,\alpha\rangle s\in\frak R,\\
\endalign
$$
and so $s'-\langle\beta,\alpha\rangle s\in S_{\sigma_{\alpha}(\beta)}$,
i.e.,
$$S_\beta-\langle\beta,\alpha\rangle S_\alpha\subset S_{\sigma_{\alpha}(\beta)}
\quad\text{for all $\alpha,\beta\in\Delta$}.
\tag S1 $$

Conversely, let $\Delta$ be a locally finite irreducible root system  in 
a vector space $V_1$ over $\Bbb Q$ with positive definite form,
and let 
$\{S_\alpha\}_{\alpha\in\Delta}$ be a family
of nonempty subsets in a vector space $V_0$ indexed by $\Delta$
satisfying (S0)
and (S1).
Extend the positive definite form  on $V_1$ to $V:=V_1\oplus V_0$
so that $V_0$ is the radical of the form.
Let 
$\frak R:=\sqcup_{\alpha\in\Delta}\ (\alpha+S_{\alpha})$.
Then  $\frak R$ satisfies the axioms (A1-4) of a LEARS.
In particular (assuming that Corollary 5 holds),  if $\Delta$ is finite and
the abelian group
$\langle\cup_{\alpha\in\Delta}\ S_{\alpha}\rangle$ is free, then
$\frak R$ is an EARS.

\proclaim{Proposition 3}
A LEARS is a directed union of EARS of finite rank.
Namely,
if $\frak R=\bigsqcup_{\alpha\in\Delta}\ (\alpha+S_{\alpha})$ 
is a LEARS
in the description above, then
$$
\frak R=\bigcup_{\Delta'}\ 
\bigcup_{\Lambda_{\Delta'}}\
\bigsqcup_{\alpha\in\Delta'}\ \big(\alpha+(\Lambda_{\Delta'}\cap S_{\alpha})\big),
$$
where $\bigcup_{\Delta'}$ means a directed union over
finite irreducible subsystems $\Delta'$ of  $\Delta$
and
$\bigcup_{\Lambda_{\Delta'}}$ means a directed union over
subgroups $\Lambda_{\Delta'}$
generated by a subset $\cup_{\alpha\in\Delta'}\ S_{\alpha}'$,
where
$S_{\alpha}'$ is chosen to be any nonempty finite subset of $S_{\alpha}$.
\endproclaim
\demo{Proof}
Note that a locally finite irreducible root system is 
a directed union of finite irreducible subsystems [LN1, Cor.3.15].
Hence $\Delta$ is a directed union of
finite irreducible subsystems $\Delta'$, and so
$\frak R$ is a directed union of $\bigsqcup_{\alpha\in\Delta'}\ (\alpha+S_{\alpha})$,
i.e.,
$\frak R=\bigcup_{\Delta'}\ 
\bigsqcup_{\alpha\in\Delta'}\ (\alpha+S_{\alpha})$.

Now, 
since $S_{\alpha}$ is a directed union of $\Lambda_{\Delta'}\cap S_{\alpha}$,
say
$S_{\alpha}=\bigcup_{\Lambda_{\Delta'}}\ 
(\Lambda_{\Delta'}\cap S_{\alpha})$,
we have
$\bigsqcup_{\alpha\in\Delta'}\ (\alpha+S_{\alpha})
=\bigcup_{\Lambda_{\Delta'}}\ \bigsqcup_{\alpha\in\Delta'}\ 
\big(\alpha+(\Lambda_{\Delta'}\cap S_{\alpha})\big)$.
Moreover,
$\langle\bigcup_{\alpha\in\Delta'}\ (\Lambda_{\Delta'}\cap S_{\alpha})\rangle
=\Lambda_{\Delta'}$
is free of finite rank.
So
$\bigsqcup_{\alpha\in\Delta'}\ 
\big(\alpha+(\Lambda_{\Delta'}\cap S_{\alpha})\big)$
is an EARS of finite rank since $\Delta'$ is a finite irreducible root system,
$\Lambda_{\Delta'}\cap S_{\alpha}$ is nonempty for all $\alpha\in\Delta'$,
and
$\{\Lambda_{\Delta'}\cap S_{\alpha}\}_{\alpha\in\Delta'}$
satisfies (S1).
\qed
\enddemo

Let us recall that we have chosen a complementary subspace $V'$ of $V^0$ to get 
$\{S_\alpha\}_{\alpha\in\Delta}$.
To classify LEARS, we now choose a nice complementary subspace.
First we define for any LEARS $\frak R$,
$$\frak R^{\red}:=
\cases
\frak R &\text{if $\frak R$ is reduced}\\
\{\alpha\in\frak R\mid\ \frac{1}{2}\alpha\notin\frak R\}
 &\text{otherwise.}
\endcases
$$

Now, note that 
$( \bar\frak R^{\red},\bar V)$ is a reduced locally finite irreducible root system.
Thus there exists a reflectable base
 $\Pi$ of $( \bar\frak R^{\red},\bar V)$
 (by Lemma 2).
We fix a preimage $\alpha\in\frak R$ for each $\bar\alpha\in\Pi$.
Let
$$
\text{$V'$ be the subspace of $V$ spanned by $\{\alpha\}_{\bar\alpha\in\Pi}$.}
$$
We call this complementary subspace a {\it reflectable subspace}
determined by 
a complete set of representatives of
a reflectable base $\Pi$ of $( \bar\frak R^{\red},\bar V)$.
Then
the subsets $S_{\alpha}$ of $V^0$ defined above satisfy the following
as in [A-P, Prop.2.11]
(see also [LN2, 4.2, 4.5, 4.10, 5.2]).

\proclaim{Lemma 4}
Let $\Delta^{\red}$ be the corresponding set to $\bar\frak R^{\red}$
determined by a reflectable subspace $V'$ as above.
Then $\Delta^{\red}\subset\frak R$, or in other words,
$$0\in S_{\alpha}
\quad\text{ for all $\alpha\in\Delta^{\red}$.}
\tag S2$$
Moreover, if $\frak R$ is reduced,
then 
$$S_{2\alpha}\cap 2S_{\alpha}=\emptyset
\quad\text{for all $2\alpha,\alpha\in\Delta$.}
\tag S3$$

\endproclaim
\demo{Proof}
For any $\alpha\in\Delta^{\red}$,
by Lemma 2  when $\bar\frak R$ is infinite,
or a well-known property for a root base when $\bar\frak R$ is finite,
$\bar\alpha=\sigma_{\bar\alpha_1}\cdots\sigma_{\bar\alpha_k}(\bar\alpha_{k+1})$
for some $\bar\alpha_1, \ldots, \bar\alpha_{k+1}\in\Pi$.
Then, $\alpha_1, \ldots, \alpha_{k+1}\in V'\cap\frak R$,
by our definition  of $V'$.
Hence 
$\sigma_{\alpha_1}\cdots\sigma_{\alpha_k}(\alpha_{k+1})\in V'$,
and
$\ol{\sigma_{\alpha_1}\cdots\sigma_{\alpha_k}(\alpha_{k+1})}
=\sigma_{\bar\alpha_1}\cdots\sigma_{\bar\alpha_k}(\bar\alpha_{k+1})
=\bar\alpha$.
So we get 
$\alpha=\sigma_{\alpha_1}\cdots\sigma_{\alpha_k}(\alpha_{k+1})$.
Therefore,
$\alpha\in\frak R$ by (A3).

For the second statement,
if $2s\in S_{2\alpha}\cap 2S_{\alpha}$ for some $s\in S_{\alpha}$,
then $2\alpha+2s\in\frak R$ and $\alpha+s\in\frak R$, contradiction.
\qed
\enddemo

Let $$G=\langle \cup_{\alpha\in\Delta}\ S_{\alpha}\rangle.$$

\proclaim{Corollary 5}
We have
$\langle \frak R\rangle=\langle \Delta\rangle\oplus G$.
In  particular,
$\langle \frak R\rangle\cap V^0=G$,
and if a LEARS $\frak R$ has nullity,
then $\langle \frak R\rangle$ is free.
\endproclaim
\demo{Proof}
Since  
$\langle \Delta\rangle=\langle \Delta^\red\rangle\subset \langle \frak R\rangle$
(by Lemma 4),
we have
$\langle \frak R\rangle=\langle \Delta\rangle\oplus G$
and
$\langle \frak R\rangle\cap V^0=G$.
Note that $\langle \Delta\rangle$ is free (by Lemma 2).
So if $G$ is free,
then $\langle \frak R\rangle$
is free.
\qed
\enddemo

Now, for a LEARS $\frak R$,
we obtain a family $\{S_\alpha\}_{\alpha\in\Delta}$
of nonempty subsets in $V^0$ satisfying
(S0), (S1) and (S2).
When $\Delta$ is a finite irreducible root system,
such a family $\{S_\alpha\}_{\alpha\in\Delta}$
satisfying (S1) and (S2)
is called a {\it root system of type $\Delta$ extended by the abelian group $G$},
and
{\it reduced} if
it satisfies (S3)
(see [Y]).

\example{Remark 6}
For $\bar\alpha\in\Pi$ and any $s_{\alpha}\in S_\alpha$,
$\alpha':=\alpha+s_{\alpha}$ is
another preimage of $\bar\alpha\in\Pi$.
Let 
$W$ be the subspace of $V$ spanned by $\{\alpha'\}_{\bar\alpha\in\Pi}$,
i.e.,
$W$ is another reflectable subspace.
Or more generally,
let $W$ be a reflectable subspace determined by a different reflectable base.
Then we get the corresponding root system $\{T_{\alpha'}\}_{\alpha'\in\Delta_W}$
extended by $G'=\langle \cup_{\alpha'\in\Delta_W}\ T_{\alpha'}\rangle$ so that
$$\frak R
=\bigsqcup_{\alpha'\in\Delta_W}\ (\alpha'+T_{\alpha'}).$$ 
The relation between $\{S_{\alpha}\}_{\alpha\in\Delta}$ and 
$\{T_{\alpha'}\}_{\alpha'\in\Delta_W}$
will be clarified in Lemma 8.
\endexample

\enskip

Root systems extended by $G$
were classified in [Y].
(The main idea comes from the classification of EARS
in [A-P].)
To explain the classification,
let us introduce some terminology.

Recall that a finite irreducible root system $\Delta$ is
one of the following types;
$\Delta=\text A_\ell $ ($\ell\geq 1$),
$\text B_\ell $ ($\ell\geq 1$, $\text B_1=\text A_1$), $\text C_\ell $ 
($\ell\geq 2$, $\text C_2=\text B_2$),
$\text D_\ell $ ($\ell\geq 4$),
$\text E_\ell $ ($\ell =6,7,8$),
$\text F_4$, $\text G_2$ or $\text{BC}_\ell $ ($\ell\geq 1$).
We partition the root system $\Delta$ according to length.
Roots of $\Delta$ of minimal length are called {\it short}.
Roots of $\Delta$ 
which are two times a short root of $\Delta$
are called {\it extra long}.
Finally, roots of $\Delta$ 
which are neither short nor extra long 
are called {\it long}.
We denote the subsets of short, long and extra long roots
of $\Delta$ by
$\Delta_{sh}$, $\Delta_{lg}$ and $\Delta_{ex}$ respectively.
Thus
$$
\Delta=\Delta_{sh}\sqcup\Delta_{lg}\sqcup\Delta_{ex}.
$$
Of course the last two terms in this union may be empty.
Indeed,
$$\Delta_{lg}=\emptyset \quad \Longleftrightarrow\quad
\text{$\Delta$ has simply laced type or type $\text{BC}_1$,}$$
and
$$\Delta_{ex}=\emptyset \quad \Longleftrightarrow\quad
\Delta=\Delta^\red.
$$
If $\Delta_{lg}\neq\emptyset$,
we use the notation $k$ for the ratio
of the long square root length to the short square root length
in $\Delta$.
Hence,
$$
k=\cases
2&\text{if 
$\Delta$ has type
$\text B_\ell $, $\text C_\ell $,
$\text F_4$ or $\text{BC}_\ell $ for $\ell\geq 2$.
}\\
3&\text{if 
$\Delta$ has type
$\text G_2$.}
\endcases
$$

For any abelian group $G$,
\roster
\item"(i)"
a subset $E$ of $G$ is called a {\it  reflection space} 
if
$E-2E\subset E$;

\item"(ii)"
a  reflection space $E$ of $G$
is called {\it full} 
if
$E$ generates $G$;

\item"(iii)"
a  reflection space $E$ of $G$ is called a {\it pointed reflection space} 
if
$0\in E$.
\endroster

These notions were introduced in [A-P] when $G$ is a full lattice
in a finite-dimensional
real vector space as a name {\it semilattice},
or earlier in a more general setting in [L].
We note that if $E$ is a full  reflection space of $G$, then
$2G+E\subset 2\langle E\rangle+E\subset E$
and so $2G+E\subset E$
(see [A-P, p.23]).
Hence, 
$$
\text{$E$ is a union of cosets of $G$ by $2G$.}
$$

Now we can state the classification of
root systems $\{S_\alpha\}_{\alpha\in\Delta}$ of type $\Delta$ extended 
by $G$
[Y, Thm 3.4]:

Set
$S_\alpha=S$ for all $\alpha\in\Delta_{sh}$,
$S_\alpha=L$ for all $\alpha\in\Delta_{lg}$
and
$S_\alpha=E$ for all $\alpha\in\Delta_{ex}$,
where $S$ is a full pointed reflection space,
$L$ is a pointed reflection space
and
$E$ is a  reflection space
satisfying
$$
\align
&L+kS\subset L,
\quad S+L\subset S,
\quad E+4S\subset E,
\\
&S+E\subset S, \quad E+2L\subset E
\quad\text{and}\quad L+E\subset L;\\
\text{moreover}, \quad& S=G \quad
\text{if $\Delta\neq \text A_1,\ \text B_\ell ,\ \text{BC}_\ell $,}\\
&\text{$L$ is a subgroup
if $\Delta=\text B_\ell \ (\ell\geq 3),\
\text F_4,\ \text G_2,\ \text{BC}_\ell \ (\ell\geq 3)$,}\\
\endalign
$$
and if $\{S_\alpha\}_{\alpha\in\Delta}$ is reduced, then
$$
E\cap 2S=\emptyset.
$$
Conversely,
let 
$S$, $L$ and $E$ be as above,
and define
$S_\alpha=S$ for all $\alpha\in\Delta_{sh}$,
$S_\alpha=L$ for all $\alpha\in\Delta_{lg}$
and
$S_\alpha=E$ for all $\alpha\in\Delta_{ex}$.
Then
$\{S_\alpha\}_{\alpha\in\Delta}$ 
is a root system extended by $G$, and if
$E\cap 2S=\emptyset$, then
$\{S_\alpha\}_{\alpha\in\Delta}$
is a reduced root system extended by $G$.
We refer to the root system $\{S_\alpha\}_{\alpha\in\Delta}$
by
$\frak R(S,L,E)_{\Delta}$.

For the case where $\Delta$ is a locally finite irreducible root system,
one can classify $\{S_\alpha\}_{\alpha\in\Delta}$
satisfying (S1) and (S2) in the same way.
In fact they were classified in [LN2, 5.9]
as extension data of locally finite root systems.
One can also obtain the classification from the fact that
$\{S_\alpha\}_{\alpha\in\Delta}$ is a directed  union
$\bigcup_{\Delta'}\ \{S_\alpha\}_{\alpha\in\Delta'}$,
where $\Delta'$ is a
finite irreducible subsystem of  $\Delta$  (see Proposition 3).
Thus the properties for $S_\alpha$ of each infinite type
$\text A_{\frak I}$,
$\text B_{\frak I}$,
$\text C_{\frak I}$,
$\text D_{\frak I}$, or
$\text{BC}_{\frak I}$
are the same as of finite type
$\text A_2$,
$\text B_3$,
$\text C_3$,
$\text D_4$, or
$\text{BC}_3$, respectively.

We note that $E\subset L\subset S$ in general,
and so $S$ spans $V^0$ by our extra condition (S0).
Moreover, from the relations
$L+kS\subset L$
and $E+4S\subset E$,
$L$ or $E$ also spans $V^0$ if it is not empty.
Thus, 
the following is known:

\proclaim{Theorem 7}
Let $\frak R$ 
be a LEARS in $V=V'\oplus V^0$
so that $\Delta$ is a locally finite irreducible root system in $V'$,
described above.
Then $E\subset L\subset S$,
$\langle S\rangle=G$,
$S$ always spans $V^0$,
and $L$ or $E$ also spans $V^0$ if it is not empty.
Moreover:

If $\Delta=\text A_{\frak I}$,
then
$\frak R=\Delta+S$,
where $S$ is a pointed reflection space of $V^0$,
and if $\text A_{\frak I}\neq \text A_1$,
then $S=G$.

If $\Delta=B_{\frak I}$,
then
$\frak R=(\Delta_{sh}+S)\sqcup (\Delta_{lg}+L)$,
where $S$ and $L$ are pointed reflection spaces of $V^0$
satisfying $2S+L\subset L$ and $S+L\subset S$,
and if $|\frak I|>2$,
then $L$ is a subgroup of $V^0$.

If $\Delta=C_{\frak I}$,
then
$\frak R=(\Delta_{sh}+S)\sqcup (\Delta_{lg}+L)$,
where $S$ and $L$ are pointed reflection spaces of $V^0$
satisfying $2S+L\subset L$ and $S+L\subset S$,
and if $|\frak I|>2$,
then $S=G$.

If $\Delta=\text D_{\frak I}$, $\text E_6$, $\text E_7$ or $\text E_8$,
then
$\frak R=\Delta+G$.

If $\Delta=\text{BC}_{\frak I}$ for $|\frak I|\geq 2$,
then
$\frak R=(\Delta_{sh}+S)\sqcup (\Delta_{lg}+L)\sqcup (\Delta_{ex}+E)$,
where $S$ and $L$ are pointed reflection spaces of $V^0$
and
$E$ is a  reflection space of $V^0$
satisfying $2S+L\subset L$, $S+L\subset S$,
$4S+E\subset E$, $S+E\subset S$,
$2L+E\subset E$ and $L+E\subset L$,
and if $|\frak I|>2$,
then $L$ is a subgroup of $V^0$.
Also,
if $\frak R$ is reduced,
then 
$E\cap 2S=\emptyset$.

If $\Delta=\text{BC}_1$,
then
$\frak R=(\Delta_{sh}+S)\sqcup (\Delta_{ex}+E)$,
where $S$ is a pointed reflection space of $V^0$
and
$E$ is a  reflection space of $V^0$
satisfying $4S+E\subset E$ and $S+E\subset S$.
Also,
if $\frak R$ is reduced,
then 
$E\cap 2S=\emptyset$.

If $\Delta=\text F_4$, then
$\frak R=(\Delta_{sh}+G)\sqcup (\Delta_{lg}+L)$,
where $L$ is a subgroup of $V^0$
satisfying $2G\subset L$.

If $\Delta=\text G_2$, then
$\frak R=(\Delta_{sh}+G)\sqcup (\Delta_{lg}+L)$,
where $L$ is a subgroup of $V^0$ satisfying $3G\subset L$.

  Conversely, each set $\frak R$ defined above is a LEARS of the specified type
  (see the paragraph right before Proposition 3).
\endproclaim

The reader should always keep in mind that
even if a LEARS $\frak R$
is reduced, the corresponding 
finite root system $\bar\frak R$ or $\Delta$ could be nonreduced.

\head
Isomorphisms
\endhead

By Theorem 7, the classification of LEARS is reduced to the classification
of triples $\{S,L,E\}$ described there.
We simply say triples, but they might be $\{S\}$, $\{S,L\}$ or $\{S,E\}$
depending on the types.
We treat these cases as special cases of triples, and we do not mention this
in the argument below.
The reader should ignore the description of $L$ or $E$ if the system does not
have $L$ or $E$,
i.e., the case $\Delta_{lg}=\emptyset$ or $\Delta_{ex}=\emptyset$.
To investigate when two triples give the same LEARS,
we show the following:
(There is a similar statement in [A1, p.577] for EARS of reduced type.)
\proclaim{Lemma 8}
In the description of Theorem 7,
let $s\in S$ and $l\in L$.
Then the triples $\{S,L,E\}$ and $\{S-s,L-l,E-2s\}$ 
give the same LEARS
(by the same $\Delta$ in Theorem 7).

Conversely, let $\{S_1,L_1,E_1\}$ be another triple obtained from a reflectable subspace $W$
of an arbitrary reflectable base.
Then,  $S_1=S-s$, $L_1=L-l$ and $E_1=E-2s$ for some $s\in S$ and $l\in L$.

\endproclaim
\demo{Proof}
Recall from the previous section that for each $\bar\alpha\in\Pi$
(a reflectable base of $\bar\frak R$),
we have considered a fixed preimage $\alpha\in\frak R$.
For each $\bar\alpha\in\Pi\cap\bar\frak R_{sh}$,
let $\alpha':=\alpha+s$,
and
for each
$\bar\alpha\in\Pi\cap\bar\frak R_{lg}$,
let $\alpha':=\alpha+l$.
Let
$U$ be the subspace of $V$ spanned by $\{\alpha'\}_{\bar\alpha\in\Pi}$.
In other words,
$U$ is another reflectable subspace.
Then 
the new family 
$\{T_{\alpha'}\}_{\alpha'\in\Delta_U}$
is a root system extended by $G$,
which gives the same LEARS.
In particular,
$\alpha+s+T_{\alpha'}=\alpha+S$
and
$\alpha+l+T_{\alpha'}=\alpha+L$.
Thus
$T_{\alpha'}=S-s$ if $\bar\alpha\in\Pi\cap\bar\frak R_{sh}$
and
$T_{\alpha'}=L-l$ if $\bar\alpha\in\Pi\cap\bar\frak R_{lg}$.
Hence, by Theorem 7, we have
$T_{\alpha'}=S-s$ for all $\alpha'\in(\Delta_U)_{sh}$
and
$T_{\alpha'}=L-l$ for $\alpha'\in(\Delta_U)_{lg}$.
Finally
(the case $\Delta_{ex}\neq\emptyset$),
for $\bar\alpha\in\Pi\cap\bar\frak R_{sh}$,
we have
$\alpha'-\alpha=s$, and so
$2\alpha'-2\alpha=2s$.
Since
$2\alpha'+T_{2\alpha'}=2\alpha+E$,
we get
$T_{2\alpha'}=2\alpha-2\alpha'+E=E-2s$.
Thus, by Theorem 7,
$T_{2\alpha'}=E-2s$
for all
$\bar\alpha,2\bar\alpha\in\bar\frak R$.

For the second statement,
let us remind the reader that
the reflectable subspace $U$ determines
 another root system $\{T_{\alpha'}\}_{\alpha'\in\Delta_U}$
extended by $G'=\langle \cup_{\alpha'\in\Delta_U}\ T_{\alpha'}\rangle$,
as in Remark 6.
Then by Theorem 7,
 the system $\{T_{\alpha'}\}_{\alpha'\in\Delta_U}$ turns out to be
 just a triple,
 that is
 $\{S_1,L_1,E_1\}$ in our assumption.
 (In particular, $G'=\langle S_1\rangle$.)
Now, for $\alpha'\in(\Delta_U)_{sh}$,
there exists $\alpha\in\Delta_{sh}$
such that $\ol{\alpha'}=\bar\alpha$.
Thus $\alpha'=\alpha+s$ for some $s\in S$.
Hence, 
$S_1=T_{\alpha'}=S-s$.
By the same argument,
we get $L_1=L-l$
for some $l\in L$.
Then, by the same argument above,
we obtain
$E'=E-2s$.
(It is enough that one of the short $\alpha'$'s satisfies $T_{\alpha'}=S-s$
and
one of the long $\alpha'$'s satisfies $T_{\alpha'}=L-l$.)
In particular, $G'=G$.
\qed
\enddemo

Two LEARS $( \frak R,V)$ and $( \frak S,W)$
are called {\it isomorphic} if there exists a linear isomorphism
$\vvp: V\longrightarrow W$ such that 
$\vvp(\frak R)=\frak S$.

\medskip

The argument to show that
 $\vvp(V^0)=W^0$ in the following lemma
 is adapted from [AY, Lemma 3.1].

\proclaim{Lemma 9}
Suppose that two LEARS are isomorphic, say
$\vvp:(\frak R,V)\tilde{\longrightarrow}(\frak S,W)$.
Then $\vvp(V^0)=W^0$ and
$\langle\vvp(\alpha),\vvp(\beta)\rangle=\langle\alpha,\beta\rangle$
for all $\alpha,\beta\in\frak R$.
Thus $\vvp$ preserves the form up to nonzero scalar.
Also, $\vvp\circ\sigma_\alpha\circ\vvp^{-1}=\sigma_{\vvp(\alpha)}$
for all $\alpha\in\frak R$.
\endproclaim
\demo{Proof}
We first show that $\vvp(V^0)=W^0$.
Let $S$ and $\Delta$ be as in Theorem 7.
Since $S$ spans $V^0$, it is enough to show that
$s\in S$ $\Rightarrow$ $s':=\vvp(s)\in W^0$,
or equivalently
 $(s',s')=0$.
Since $S\pm 2S\subset S$, we have $ns\in S$ for all $n\in\Bbb Z$.
Let $\alpha\in\Delta_{sh}$ and $\alpha':=\vvp(\alpha)$.
By Theorem 7,
$\alpha+ns\in\frak R$ for all $n\in\Bbb Z$
and so 
$\alpha'+ns'=\vvp(\alpha+ns)\in\vvp(\frak R)=\frak S$.
But then by the axiom (A2) of the definition of a LEARS,
we have
$$
\langle \alpha',\alpha'+ns'\rangle
=\frac{2(\alpha',\alpha'+ns')}{(\alpha'+ns',\alpha'+ns')}
=\frac{2(\alpha',\alpha')+2n(\alpha',s')}
{(\alpha',\alpha')+2n(\alpha',s')
+n^2(s',s')}\in\Bbb Z
$$
for all $n\in\Bbb Z$
which implies
$(s',s')=0$
(note that  $(\alpha',\alpha')\neq 0$ and let $n\rightarrow\infty$).
Thus we have shown that 
$\vvp(V^0)=W^0$.
Then $\vvp$ induces a linear isomorphism
$\bar\vvp:\bar V\tilde{\longrightarrow}\bar W$
with $\bar\vvp(\bar\frak R)= \bar\frak S$,
and this is what means an isomorphism of
locally finite root systems in [LN1].
Thus, by [LN1, Lem.3.7],
we have
$\langle\bar\vvp(\bar\alpha),\bar\vvp(\bar\beta)\rangle=\langle\bar\alpha,\bar\beta\rangle$
for all $\bar\alpha,\bar\beta\in\bar\frak R$.
So the second statement is shown since we always have
$(\alpha,\beta)=(\bar\alpha,\bar\beta)$.
The third statement follows from the equivalence between
connectedness and irreducibility in our systems (see [LN2, Lem.2.7]
or [MP, Prop.3.4.6]).
The last statement is now clear.
\qed
\enddemo

We introduce a notion of similarity for triples
following [A-P].

Let  $(S_1,L_1,E_1)$ and $(S_2,L_2,E_2)$
be two triples satisfying the properties in Theorem 7
in vector spaces $W_1$ and $W_2$, respectively.
We say 
that $(S_1,L_1,E_1)$ and $(S_2,L_2,E_2)$
are {\it similar}, denoted
$(S_1,L_1,E_1)\sim(S_2,L_2,E_2)$,
if
there exists an isomorphism $\vvp$ from $W_1$ onto $W_2$
such that
$\vvp(S_1)=S_2-s_2$,
$\vvp(L_1)=L_2-l_2$ and 
$\vvp(E_1)=E_2-2s_2$
for some $s_2\in S_2$ and $l_2\in L_2$.
The similarity is an equivalence relation.

The following theorem says that
there is a 1-1 correspondence between the isomorphism classes of LEARS
and the similarity classes of triples.
The theorem generalizes [A-P, Thm 3.1] and our proof is simpler.
\proclaim{Theorem 10}
Suppose that
 $\vvp:(\frak R_1,V_1;V_1',\Delta_1;S_1,L_1,E_1)\tilde{\longrightarrow}
(\frak R_2,V_2;V_2',\Delta_2;S_2,L_2,E_2)$
is an isomorphism of LEARS.
Let 
$$
\text{$\zeta:=$(projection onto $V_2'$) $\circ$ $\vvp\mid_{V_1'}$
\ and\ \
$\psi:=$(projection onto $V_2^0$) $\circ$ $\vvp\mid_{V_1'}$.}
$$
Then
$\zeta:(\Delta_1,V_1')\tilde{\longrightarrow}(\Delta_2,V_2')$
and $\vvp(S_1)=S_2-s_2$,  $\vvp(L_1)=L_2-l_2$
and  $\vvp(E_1)=E_2-2s_2$ for some
$s_2\in S_2$ and $l_2\in L_2$.

Conversely,
if $\zeta:(\Delta_1,V_1)\tilde{\longrightarrow}(\Delta_2,V_2)$
is an isomorphism of locally finite irreducible root systems,
two triples $(S_1,L_1,E_1)$ in a vector space $W_1$
and $(S_2,L_2,E_2)$ in a vector space $W_2$
satisfy the conditions in Theorem 7 depending on the type of $\Delta_1$,
and 
$\vvp$ is an isomorphism from $W_1$ onto $W_2$
so that
$\vvp(S_1)=S_2-s_2$,  $\vvp(L_1)=L_2-l_2$
and  $\vvp(E_1)=E_2-2s_2$ for some
$s_2\in S_2$ and $l_2\in L_2$,
then
$(\frak R(S_1,L_1,E_1),V_1\oplus W_1)$
is isomorphic to
$(\frak R(S_2,L_2,E_2),V_2\oplus W_2)$.
\endproclaim
\demo{Proof}
We have $\zeta(\Delta_1), \Delta_2\subset V_2'$
and so $\ol{\zeta(\Delta_1)}=\bar\frak R_2=\bar\Delta_2$
since $\vvp(V_1^0)=V_2^0$ (Lemma 9).
Hence
$\zeta(\Delta_1)=\Delta_2$,
and so $\zeta$ is an isomorphism of the root systems.
Also, for a fixed $\alpha\in(\Delta_1)_{sh}$,
we have
$\vvp(\alpha+S_1)=\zeta(\alpha)+\psi(\alpha)+\vvp(S_1)\subset\frak R_2$,
and so
$\psi(\alpha)+\vvp(S_1)=S_2$
since $\zeta(\alpha)\in(\Delta_2)_{sh}$.
Also,  $s_2:=\psi(\alpha)\in S_2$ since $0\in\vvp(S_1)$.
Similarly, for a fixed $\beta\in(\Delta_1)_{lg}$,
we get
$l_2+\vvp(L_1)=L_2$
for  $l_2:=\psi(\beta)\in L_2$.
Finally, if $2\alpha\in\Delta_1$,
then
$\vvp(2\alpha+E_1)=\zeta(2\alpha)+2\psi(\alpha)+\vvp(E_1)$,
and so
$2s_2+\vvp(E_1)=E_2$.

For the second statement,
let $\tilde\vvp=\zeta\oplus\vvp$.
Then 
$$
\align
\tilde\vvp:(\frak R(S_1,L_1,E_1),V_1\oplus W_1)
&\tilde{\longrightarrow}
(\frak R(\vvp(S_1),\vvp(L_1),\vvp(E_1)),V_2\oplus W_2)\\
&=(\frak R(S_2-s_2,L_2-l_2,E_2-2s_2),V_2\oplus W_2)\\
&=(\frak R(S_2,L_2,E_2),V_2\oplus W_2)
\quad\text{by Lemma 8}. \qed\\
\endalign
$$

\enddemo

\example{Remark 11}
If two LEARS
 $\frak R_1$ and $\frak R_2$
are isomorphic,
then
$\langle S_1\rangle/\langle L_1\rangle$
and $\langle S_2\rangle/\langle L_2\rangle$
are clearly isomorphic as abelian groups.
Also the reducibility of LEARS is an isomorphic invariant.

\endexample

\medskip

\head
Special case
\endhead

We consider LEARS of null dimension 1.
Then the abelian group $G$ in Theorem 7 is just a subgroup of $\Bbb Q$.
We first observe special properties for
a cyclic group or a subgroup of $\Bbb Q$.
Let us recall the concept of divisibility for an arbitrary abelian group $G$.
We say that a prime number {\it $p$ is divisible in $G$}
or {\it $G$ is divisible by $p$}
if $G=pG$, or equivalently $px=g$ has a solution $x$ in $G$ for any $g\in G$.
Any cyclic group of infinite order is not divisible by any prime.
The following is a useful exercise ([G, p.8]):
$$
\text{If $mx=ng$ for $(m,n)=1$ has a solution $x$ in $G$,
then $my=g$ has a solution $y$ in $G$.}
\tag$*$
$$

\proclaim{Lemma 12}
(1)
If $S$ is a full  reflection space of a cyclic group $G$,
then $S=G$ or $S=2G+s$ for any $s\in G\setminus 2G$.
So if $S$ is a full pointed reflection space of a cyclic group $G$,
then $S=G$.

(2)
Suppose that $G$ is a subgroup of $\Bbb Q$.
If $G$ is
not divisible by a prime $p$,
then $G/p^nG\cong\Bbb Z_{p^n}$ for any $n\in\Bbb N$.
Moreover, if $G/H\cong\Bbb Z_{p^n}$ for some subgroup $H$ of $G$
and some $n\in\Bbb N$,
then $G$ is not divisible by $p$
and $H=p^nG$.

(3)
If $S$ is a full  reflection space of a subgroup $G$ of $\Bbb Q$
divisible by $2$,
then $S=G$.

(4)
The same statement in (1) is true for 
a subgroup $G$ of $\Bbb Q$ not divisible by $2$.

\endproclaim
\demo{Proof}
For (1),
we have $G=2G\sqcup (2G+s)$ for any $s\in G\setminus 2G$
if $G\neq 2G$.
(Note that 
$G$ is finite of odd order
$\Leftrightarrow$ $G= 2G$.)
Since $S$ is full,
$S$ is a union of cosets of $G$ by $2G$,
and $S\neq 2G$ if $G\neq 2G$.
So $S=G$ or $S=2G+s$.
Moreover, $2G\subset S$ if $0\in S$, and hence (1) is proved.

For (2),
by the divisibility,
there exists $g\in G\setminus pG$.
We claim that
$0, g,2g, \ldots,(p^n-1)g$ are distinct
modulo $p^nG$.
($G$ can be any torsion free group for this claim.)
Suppose that two of them are equal.
Then $p^rqg=p^ng'$ for some $r<n$, $(p,q)=1$ and $g'\in G$.
Hence $p^{n-r}g'=qg$ (since $G$ is torsion free).
Then by ($*$) above, $py=g$ has a solution $y$ in $G$,
which contradicts our choice of $g$.
Thus we showed the claim, and
the order of $g$ in $G/p^nG$ is $p^n$.

Let $g'\in G$.
Since $G':=\langle g,g'\rangle$ is cyclic,
we have  $G'/p^nG'\cong \Bbb Z_{p^n}$.
Hence,
$G'/p^nG'=\langle g\rangle$, and 
$g'$ is equal to one of $g,2g, \ldots,(p^n-1)g$ or $p^ng$ modulo $p^nG'$,
and so is in the modulo $p^nG$
since $p^nG'\subset p^nG$.
Hence
$G/p^nG=\{0, g,2g, \ldots,(p^n-1)g\}$,
which is a cyclic group with $p^ng=0$.
But since the order of $g$ in $G/p^nG$ is $p^n$,
we obtain
$G/p^nG\cong \Bbb Z_{p^n}$.

For the second statement,
if $G$ is divisible by $p$,
then for any $g\in G$, $g=p^ng'$ for some $g'\in G$.
So for any $\bar g\in G/H$, $\bar g=\ol{p^ng'}=\bar 0$, 
which means $G/H=0$, contradiction.
Hence, $G$ is not divisible by $p$.
Thus by the first statement, we have $G/p^nG\cong \Bbb Z_{p^n}$.
Note that $G/H\cong \Bbb Z_{p^n}$ implies $p^nG\subset H$.
So there is a natural epimorphism $\pi$ from
$G/p^nG$ onto $G/H$.
But the order of both groups is $p^n$,
and hence $\pi$ is an isomorphism and $p^nG= H$.

For (3), 
we have 
$G=G+S=2G+S\subset S$, and hence $G=S$.

For (4),
applying (2) for $p=2$,
we have $G=2G\sqcup (2G+s)$ for any $s\in G\setminus 2G$.
Thus we are done.
(We note that 
since $s\in\langle 2G+s\rangle$,
we have $2G\subset\langle 2G+s\rangle$,
and hence $\langle 2G+s\rangle=G$.)
\qed
\enddemo

We will use the special cases of Lemma 12(2) later,
namely $p=2$ or $p=3$ for $n=1$.
Note that this is a special property of subgroups of $\Bbb Q$.
For example,
if $G=\langle\sqrt 2,\sqrt 3\rangle\subset \Bbb R$,
then
$G/2G\cong\Bbb Z_2\times\Bbb Z_2$.

\proclaim{Corollary 13}
Let $\frak R$ 
be a LEARS of null dimension $1$
in $V=V'\oplus V^0$,
$\Delta$ a locally finite irreducible root system in $V'$,
and $G$ a subgroup in $V^0=\Bbb Q$,
as described above.

(1a)
For the case where $\frac{1}{2}\notin G$:

If $\Delta=\text A_{\frak I}$, $\text D_{\frak I}$,
$\text E_6$, $\text E_7$ or $\text E_8$
then
$\frak R=\Delta+ G$.

If $\Delta=B_{\frak I}$, $C_{\frak I}$
or $\text F_4$,
then $\frak R=\Delta+ G$ or
$\frak R=(\Delta_{sh}+ G)\sqcup (\Delta_{lg}+2 G)$.

If $\Delta=\text G_2$,
then $\frak R=\Delta+ G$ or
$\frak R=(\Delta_{sh}+ G)\sqcup (\Delta_{lg}+3 G)$.

If $\Delta=\text{BC}_{\frak I}$ for $|\frak I|>1$,
then
$$
\align
\frak R&=\Delta+ G,\\
\frak R&=((\Delta_{sh}\sqcup \Delta_{lg})+ G)\sqcup (\Delta_{ex}+2 G),\\
\frak R&=(\Delta_{sh}+ G)\sqcup ((\Delta_{lg}\sqcup \Delta_{ex})+2 G),\\
\frak R&=(\Delta_{sh}+ G)\sqcup (\Delta_{lg}+2 G)\sqcup (\Delta_{ex}+4 G)
\quad\text{or}\\
\frak R&=((\Delta_{sh}\sqcup \Delta_{lg})+ G)\sqcup (\Delta_{ex}+2 G+s)
\quad\text{for any $s\in G\setminus 2G$,}
\endalign
$$
and moreover,
if $\frak R$ is reduced,
then 
only the last case happens.

If $\Delta=\text{BC}_1$,
then
$$
\align
\frak R&=\Delta+ G,\\
\frak R&=(\Delta_{sh}+G)\sqcup (\Delta_{ex}+2 G),\\
\frak R&=(\Delta_{sh}+ G)\sqcup (\Delta_{ex}+4 G)
\quad\text{or}\\
\frak R&=(\Delta_{sh}+G)\sqcup (\Delta_{ex}+2 G+s)
\quad\text{for any $s\in G\setminus 2G$,}
\endalign
$$
and moreover,
if $\frak R$ is reduced,
then 
only the last case happens.

(1b)
If $G$ is divisible by $2$,
then $\frak R=\Delta+ G$ in any type of $\Delta$.

(2) If $\frak R$ is a LARS, then
$\frak R$ has
the same description as in (1a)
by changing $G$ into $\Bbb Zs$,
where $s\in G$ so that $G=\Bbb Zs$.
\endproclaim

\demo{Proof}
First of all, note that all the LEARS in the list above are not isomorphic
by Remark 11.
Also, by Lemma 12,
we always have $G=S$ and
$L$ is a group since $S$ is a full pointed reflection space of $G\subset\Bbb Q$
and $L$ is a full pointed reflection space of $\langle L\rangle\subset\Bbb Q$, and
hence $2G\subset L$ (or $3G\subset L$ for type $\text G_2$)
by Theorem 7.
But then by Lemma 12,
$L=2G$ or $G$ ($L=3G$ or $G$ for type $\text G_2$).
So we are done except for the type $\text{BC}_{\frak I}$.

Now for $|\frak I|>1$, if $L=G$, then  
$2G+E\subset E$, and so $E$ is a union of cosets of $G$ by $2G$.
Hence, 
$E=2G$, $G$ or $2G+s$ for any $s\in G\setminus 2G$.

If $L=2G$, then 
$E\subset 2G$.
So we have $4G+E\subset E\subset 2G$,
and hence $E$ is a union of cosets of $2G$ by $4G$.
Hence, $E=4G$, $2G$ or $4G+g$ for any $g\in 2G\setminus 4G$,
and $4G+g=4G+2s$ for any $s\in G\setminus 2G$.
But  $4G+2s$ is excluded since
$(G,2G,4G)$ and $(G,2G,4G+2s)$ are similar.
Also,
$E=2G+s$ is the only reduced one
since others do not satisfy $E\cap 2S=\emptyset$.

For the type $\text{BC}_1$,
we only have $4G+E\subset E$,
and so $E$ is a union of cosets of $G$ by $4G$.
By Lemma 12, we have 
$G=4G\sqcup (4G+s)\sqcup (4G+2s)\sqcup(4G+3s)$
for any $s\in G\setminus 2G$.
Note that $4G+s\subset E$ $\Leftrightarrow$ $4G+3s\subset E$
since $2E+E\subset E$.
Also, $4G, 4G+s\subset E$  $\Rightarrow$ $E=G$.
Hence,
$E=4G$, $4G+s$, $4G+2s$,
$4G\sqcup (4G+2s)=2G$ or $G$.
But if $E=4G+s$,
then $0\notin E$ and $\langle E\rangle=G$,
and so $E=2G+s$.
As above,
$4G+2s$ is excluded since
$(G,4G)$ and $(G,4G+2s)$ are similar,
and
$E=2G+s$ is the only reduced one
since others do not satisfy $E\cap 2S=\emptyset$.

(1b):
We have $S=G$ by Lemma 12(3).
Moreover, by Theorem 7,
we have $L\supset 2S+L=2G+L=G+L=G$ and 
$E\supset 2L+E=2G+E=G+E=G$,
and hence $G=S=L=E$
(if $L$ or $E$ is empty).
This shows (1b).

For (2),
we have $\langle \frak R\rangle\cap V^0
=\langle S\rangle$
has rank 1, and so
there exists $s\in S$ so that $S=\Bbb Zs=G$
(see Lemma 12(1)).
\qed
\enddemo

\example{Remark 14}
(1)  Nonreduced EARS of nullity 1, 2 and 3 were already classified
in [AKY].

(2)  
Note that a free abelian group is not divisible by any $p$.
An example of a subgroup of $\Bbb Q$ not divisible by $p$,
which is not free,
is the localization $\Bbb Z_{(p)}$ of $\Bbb Z$ 
by the prime ideal $(p)=p\Bbb Z$.
Also, $\Bbb Z[\frac{1}{q}]=\langle\frac{1}{q^n}\mid n\in\Bbb N\rangle$
for any prime $q$ different from $p$
is another example of a subgroup of $\Bbb Q$ not divisible by $p$,
which is not free.
Note that $\Bbb Z[\frac{1}{q}]\subset \Bbb Z_{(p)}$
and that $\Bbb Z_{(p)}$ and $\Bbb Z[\frac{1}{q}]$
are not just subgroups but subrings of $\Bbb Q$.
There are some examples which are not subrings.
For example,
$\Bbb Z_{(p)}+\langle\frac{1}{p^n}\rangle$
is neither divisible by $p$
nor a subring of $\Bbb Q$ (nor free).
Note that $\Bbb Z_{(p)}\subset\Bbb Z_{(p)}+\langle\frac{1}{p}\rangle
\subset\Bbb Z_{(p)}+\langle\frac{1}{p^2}\rangle\subset \cdots$.
Also,
$\langle\frac{1}{p_1},\frac{1}{p_2},\ldots\rangle$
for any infinite series of distinct primes $p_1, p_2,\ldots$
is an example of a subgroup of $\Bbb Q$ not divisible by $p$ 
and not a subring of $\Bbb Q$ (and not free, even if one of the $p_i$'s is equal to $p$).
Note that the torsion-free abelian groups of rank 1 were classified (but not for rank $>1$).
\endexample

\medskip

We note that there are 14 reduced irreducible affine root systems,
i.e.,
$\text A_\ell ^{(1)}$,
$\text B_\ell ^{(1)}$, 
$\text B_\ell ^{(2)}$, 
$\text C_\ell ^{(1)}$, 
$\text C_\ell ^{(2)}$,
$\text D_\ell ^{(1)}$,
$\text{BC}_\ell ^{(2)}$,
$\text E_6^{(1)}$,
$\text E_7^{(1)}$,
$\text E_8^{(1)}$,
$\text F_4^{(1)}$, 
$\text F_4^{(2)}$, 
$\text G_2^{(1)}$
and
$\text G_2^{(3)}$, 
 by Moody's Label,
 and correspondingly there are 14 affine Lie algebras.
It is worth mentioning that there are 14 reduced LARS
from Corollary 13,
and they are obtained by just changing $\ell$ of the first 7 above into an infinite index set
$\frak I$.
For the convenience of the reader, we summarize this remark with the above label,
denoting the specific type instead of $\Delta$ and identifying $\Bbb Zs$ with $\Bbb Z$
in Corollary 13:

\proclaim{Corollary 15}
There are only seven reduced LARS of infinite rank.
Namely, 
$$
\align
\text A_{\frak I}^{(1)}
&=\text A_{\frak I}+\Bbb Z,\\
\text B_{\frak I}^{(1)}
&=\text B_{\frak I}+\Bbb Z,\\
\text C_{\frak I}^{(1)}
&=\text C_{\frak I}+\Bbb Z,\\
\text D_{\frak I}^{(1)}
&=\text D_{\frak I}+\Bbb Z,\\
\text B_{\frak I}^{(2)}
&=\big((\text B_{\frak I})_{sh}+ \Bbb Z\big)\sqcup \big((\text B_{\frak I})_{lg}+2 \Bbb Z\big),\\
\text C_{\frak I}^{(2)}
&=\big((\text C_{\frak I})_{sh}+ \Bbb Z\big)\sqcup \big((\text C_{\frak I})_{lg}+2 \Bbb Z\big)
\quad\text{and}\\
\text{BC}_{\frak I}^{(2)}
&=\bigg(\big(({BC}_{\frak I})_{sh}\sqcup({BC}_{\frak I})_{lg}\big)+ \Bbb Z\bigg)
\sqcup \big(({BC}_{\frak I})_{ex}+2 \Bbb Z+1\big).\\
\endalign
$$
\endproclaim

\head
Locally $(G, \tau)$-loop algebras
\endhead

We give examples of Lie algebras whose root systems are LEARS of null dimension 1.
All algebras and tensors are over a field $F$ of characteristic 0.
Let ${\frak I}$ be any index set. The locally finite split simple Lie algebra of type 
$\text X_{\frak I}$ 
(introduced in [NS])
is defined as a subalgebra of the matrix algebra $gl_{\frak I}(F)$,
$gl_{2{\frak I}+1}(F)$ or $gl_{2{\frak I}}(F)$ consisting of matrices having only a finite number of nonzero
entries:
(There is a more general construction in [N].)

Type $\text A_{\frak I}$;
$sl_{\frak I}(F)=\{x\in gl_{\frak I}(F)\mid \tr (x)=0\}$;

Type $\text B_{\frak I}$;
$o_{2{\frak I}+1}(F)=\{x\in gl_{2{\frak I}+1}(F)\mid sx=-x^ts\}$;

Type $\text C_{\frak I}$;
$sp_{2{\frak I}}(F)=\{x\in gl_{2{\frak I}}(F)\mid s_-x=-x^ts_-\}$;

Type $\text D_{\frak I}$;
$o_{2{\frak I}}(F)=\{x\in gl_{2{\frak I}}(F)\mid s_+x=-x^ts_+\}$;

\noindent where 
 $x^t$ is the transpose of $x$,
$$s=
\left (\matrix \format\r&\quad\r &\quad\r    \\
   0   &I_{\frak I}    & 0   \\
  I_{\frak I}   & 0    & 0  \\
   0   & 0   &   1   
\endmatrix
\right ),\
s_-=\left (\matrix \format\r&\quad\r &\quad\r &\quad\r   \\
   0   &I_{\frak I}       \\
  -I_{\frak I}   & 0   \\
\endmatrix
\right )
\quad\text{or}\quad
s_+=\left (\matrix \format\r&\quad\r &\quad\r &\quad\r   \\
   0   &I_{\frak I}       \\
  I_{\frak I}   & 0   \\
\endmatrix
\right ),
$$
and $I_{\frak I}$ is the identity matrix of size ${\frak I}$.
(Each Lie algebra of type $\text X_{\frak I}$ has the locally finite irreducible root system 
of type $\text X_{\frak I}$ [NS].)

Let $G=(G,+,0)$ be an abelian group.
Let 
$$F^\tau[G]=F^\tau[G,t]=\bigoplus_{g\in G}\ Ft^g$$ 
be a twisted commutative group algebra of $G$
with symmetric twisting $\tau:G\times G\longrightarrow F^\times$, i.e.,
$$\tau(g,h)=\tau(h,g)
\quad\text{and}\quad 
\tau(g+h,k)\tau(g,h)=\tau(g,h+k)\tau(h,k)$$ 
so that
$$t^gt^h=\tau(g,h)t^{g+h}$$ for all $g,h,k\in G$.
We call the following four Lie algebras {\it locally untwisted $(G, \tau)$-loop algebras},
and {\it untwisted $(G, \tau)$-loop algebras} if $\frak I$ is finite.

Type $\text A_{\frak I}^{(1)}$;
$sl_{\frak I}(F)\otimes F^\tau[G]$;

Type $\text B_{\frak I}^{(1)}$;
$o_{2{\frak I}+1}(F)\otimes F^\tau[G]$;

Type $\text C_{\frak I}^{(1)}$;
$sp_{2{\frak I}}(F)\otimes F^\tau[G]$;

Type $\text D_{\frak I}^{(1)}$;
$o_{2{\frak I}}(F)\otimes F^\tau[G]$.

Also, for each finite-dimensional split
simple Lie algebra  $\frak g$ over $F$
of type $\text E_6$, $\text E_7$, $\text E_8$, $\text F_4$ or $\text G_2$,
we call the Lie algebra
$\frak g\otimes F^\tau[G]$
an {\it untwisted $(G, \tau)$-loop algebra} of type
$\text E_6^{(1)}$, 
$\text E_7^{(1)}$, 
$\text E_8^{(1)}$, 
$\text F_4^{(1)}$ or 
$\text G_2^{(1)}$.

If there exists a subgroup $G'$ so that $G/G'\cong \Bbb Z_2$,
then $G=G'\sqcup (G'+g_1)$ for any $g_1\in G\setminus G'$,
and so
$F^\tau[G]=F^\tau[G']\oplus t^{g_1}F^\tau[G']$.
(For example, take any subgroup $G$ of $\Bbb Q$ which is not divisible by 2,
and $G':=2G$, by Lemma 12(2).)
In this case we call the following three Lie algebras
{\it locally twisted $(G, \tau)$-loop algebras},
and {\it twisted $(G, \tau)$-loop algebras} if $\frak I$ is finite.
(There is a way to construct by Kac, using an automorphism 
of a Lie algebra in [K, Ch.8].
But we chose the following way
by [BZ] and [ABG] since this construction can be generalized to nonassociative
coordinates and is simpler.)

(1) Type $\text B_{\frak I}^{(2)}$;
$(o_{2{\frak I}+1}(F)\otimes F^\tau[G'])\oplus (V\otimes  t^{g_1}F^\tau[G'])$,

\noindent
where $V=F^{(2{\frak I}+1)}$ is the natural $o_{2{\frak I}+1}(F)$-module;

(2) Type $\text C_{\frak I}^{(2)}$;
$(sp_{2{\frak I}}(F)\otimes F^\tau[G'])\oplus (\frak s_-\otimes  t^{g_1}F^\tau[G'])$,

\noindent
where $\frak s_-=\{x\in sl_{2{\frak I}}(F)\mid s_-x=x^ts_-\}$;

(3) Type $\text {BC}_{\frak I}$;
$(o_{2{\frak I}+1}(F)\otimes F^\tau[G'])\oplus (\frak s\otimes  t^{g_1}F^\tau[G'])$,

\noindent
where $\frak s=\{x\in sl_{2{\frak I}+1}(F)\mid sx=x^ts\}$.

Note that $sl_{2{\frak I}}(F)=sp_{2{\frak I}}(F)\oplus \frak s_-$
and
$sl_{2{\frak I}+1}(F)=o_{2{\frak I}+1}(F)\oplus \frak s$.

The Lie bracket of each untwisted type is natural,
i.e., $[x\otimes t^g, y\otimes t^h]=[x,y]\otimes \tau(g,h)t^{g+h}$.
The Lie bracket of type $\text C_{\frak I}^{(2)}$ or $\text {BC}_{\frak I}$ is also natural
since 
$$
\align
&[sp_{2{\frak I}}(F), \frak s_-]\subset  \frak s_-,
\quad
[\frak s_-,\frak s_-]\subset sp_{2{\frak I}}(F),\\
&[o_{2{\frak I}+1}(F), \frak s]\subset  \frak s
\quad\text{and}\quad
[\frak s,\frak s]\subset o_{2{\frak I}+1}(F).
\endalign
$$
Note that $\text C_{\frak I}^{(2)}$ is a subalgebra of $sl_{2{\frak I}}(F)\otimes F^\tau[G]$,
and $\text {BC}_{\frak I}$ is a subalgebra of $sl_{2{\frak I}+1}(F)\otimes F^\tau[G]$.

For $\text B_{\frak I}^{(2)}$,
we have $o_{2{\frak I}+1}(F)V\subset V$, and so we define
the bracket of $o_{2{\frak I}+1}(F)$ and $V$ by the natural action,
i.e.,
$[x,v]=xv=-[v,x]$ for $x\in o_{2{\frak I}+1}(F)$ and $v\in V$.
However, there is no bracket on $V$.
So we define a bracket on $V$
so that $[V,V]\subset o_{2{\frak I}+1}(F)$ as follows.
First,
let $(\cdot,\cdot)$ be the bilinear form on $V$ determined by $s$.
Then there is a natural identification
$$
o_{2{\frak I}+1}(F)=D_{V,V}:=\spa_F\{D_{v,v'}\mid v,v'\in V\},
$$
where 
$D_{v,v'}\in\End (V)$ is defined by
$D_{v,v'}(v'')=(v',v'')v-(v,v'')v'$ for $v''\in V$.
Thus we define
$[v,v']:=D_{v,v'}$.
Note that $[v',v]=-[v,v']$.
It is easy to check that the bracket
$$
\align
&[x\otimes t^{g}+v\otimes  t^{g'+g_1},x'\otimes t^{h}+v'\otimes  t^{h'+g_1}]\\
=&[x,x']\otimes \tau(g,h)t^{g+h}+D_{v,v'}\otimes \tau(g'+g_1,h'+g_1)t^{g'+h'+2g_1}\\
&+xv'\otimes \tau(g,h'+g_1)t^{g+h'+g_1}-x'v\otimes \tau(g'+g_1,h)t^{g'+h+g_1}
\endalign
$$
defines a Lie bracket for $g,g',h,h'\in G'$,
$x,x'\in o_{2{\frak I}+1}(F)$,
$v,v'\in V$.

Also, we define two more twisted $(G, \tau)$-loop algebras.
(We use the way by Kac [K, Ch.8] for  $\text F_4^{(2)}$
in order to avoid introducing a 27-dimensional exceptional Jordan algebra.
But for $\text G_2^{(3)}$, we again use the way in [BZ]
since we do not need to assume the existence of a primitive cubic root of unity
in our base field $F$.)

(4) Type $\text F_4^{(2)}$:
Assume that $F^\tau[G]=F^\tau[G,t]=F^\tau[G']\oplus t^{g_1}F^\tau[G']$ 
with $g_1\in G\setminus G'$ again.
Let $\frak g$ be the finite-dimensional split simple Lie algebra of type $\text E_6$,
and $\sigma$ be the automorphism of $\frak g$ of order 2
determined by the diagram automorphism.
Define the automorphism $\tilde\sigma$ of $\text E_6^{(1)}=\frak g\otimes F^\tau[G]$
by $\tilde\sigma(x\otimes t^{g_1})=-\sigma(x)\otimes t^{g_1}$.
The subalgebra $L(\text F_4^{(2)},F^\tau[G])$ of $\text E_6^{(1)}$ fixed by $\tilde\sigma$
is called a {\it twisted $(G, \tau)$-loop algebra} of type $\text F_4^{(2)}$.
We note that the subalgebra $\frak g'$ of $\frak g$ fixed by $\sigma$ has type $\text F_4$,
say
$\frak g'=\bigoplus_{\mu\in\text F_4\cup\{0\}}\ \frak g'_\mu$.
Let $\frak s$ be the $(-1)$-eigenspace.
Then $\frak s$ is an irreducible highest weight $\frak g'$-module
whose highest weight is the highest short root in $\text F_4$.
Thus $L(\text F_4^{(2)},F^\tau[G,t])
=\bigoplus_{\mu\in\text F_4\cup\{0\}}\ (\frak g'_\mu\otimes F^\tau[G'])
\oplus \bigoplus_{\mu\in(\text F_4)_{sh}}\ (\frak s_\mu\otimes t^{g_1}F^\tau[G'])$.

(5) Type $\text G_2^{(3)}$:
Assume this time that 
 there exists a subgroup $G'$ so that $G/G'\cong \Bbb Z_3$.
Then $G=G'\sqcup (G'+g_1)\sqcup (G'+2g_1)$ for any $g_1\in G\setminus G'$,
and so
$$F^\tau[G]=F^\tau[G']\oplus t^{g_1}F^\tau[G']\oplus t^{2g_1}F^\tau[G'].$$
(For example, take any subgroup $G$ of $\Bbb Q$ which is not divisible by 3,
and $G':=3G$, by Lemma 12(2).)
As in [ABGP] (or in [BZ]),
let $\Bbb O$ be a split octonion algebra over $F$,
and $t:\Bbb O\longrightarrow F$ the normalized trace on $\Bbb O$,
in which $\Bbb O=F1\oplus\Bbb O_0$, where
$\Bbb O_0=\{x\in\Bbb O\mid t(x)=0\}$.
Moreover, if $x,y\in\Bbb O$, we have
$xy=t(xy)1+x*y$
for some unique $x*y\in\Bbb O_0$.
One can check that $x*y=-y*x$ for $x,y\in\Bbb O_0$.
Next, let 
$$
D_{\Bbb O,\Bbb O}:=\spa_F\{D_{x,y}\mid x,y\in\Bbb O\},
$$
where 
$D_{x,y}=\frac{1}{4}(L_{[x,y]}-R_{[x,y]}-3[L_x,R_y])$.
(Here $L_x$ and $R_x$ denote the left and right multiplication operators 
by $x$ in $\Bbb O$.)
Then
$D_{\Bbb O,\Bbb O}$
is the Lie algebra of all derivations of $\Bbb O$ and
$D_{\Bbb O,\Bbb O}$
is a split simple Lie algebra of type $\text G_2$ over $F$.
Let
$$
L(\text G_2^{(3)},F^\tau[G,t])=(D_{\Bbb O,\Bbb O}\otimes F^\tau[G'])
\oplus (\Bbb O_0\otimes  t^{g_1}F^\tau[G'])
\oplus (\Bbb O_0\otimes  t^{2g_1}F^\tau[G']).
$$
One can check that the bracket
$$
\align
&[D\otimes t^{g}+x\otimes  t^{g'+g_1}+x'\otimes  t^{g''+2g_1},
D'\otimes t^{h}+y\otimes  t^{h'+g_1}+y'\otimes  t^{h''+2g_1}]\\
&=[D,D']\otimes \tau(g,h)t^{g+h}
+Dy\otimes \tau(g,h'+g_1)t^{g+h'+g_1}
+Dy'\otimes \tau(g,h''+2g_1)t^{g+h''+2g_1}\\
&-D'x\otimes \tau(g'+g_1,h)t^{g+h+g_1}
+(x*y)\otimes \tau(g'+g_1,h'+g_1)t^{g'+h'+2g_1}\\
&+D_{x,y'}\otimes \tau(g'+g_1,h''+2g_1)t^{g'+h''+3g_1}
-D'x'\otimes \tau(g''+2g_1,h)t^{g''+h+2g_1}\\
&+D_{x',y}\otimes \tau(g''+2g_1,h'+g_1)t^{g''+h'+3g_1}
+(x'*y')\otimes \tau(g''+2g_1,h''+2g_1)t^{g''+h''+4g_1}
\endalign
$$
defines a Lie bracket for 
$D,D'\in D_{\Bbb O,\Bbb O}$,
$x,x',y,y'\in \Bbb O_0$ and $g,g',g'', h,h', h''\in G'$.
In fact, if we define an $F$-linear map $\tr$ on $F^\tau[G,t]$
by
$$
\tr(t^g)=\cases
t^g&\text{if 
$g\in G'$}\\
0&\text{otherwise}
\endcases
$$
(so $\tr$ is
 an $F^\tau[G']$-linear map on the algebra $F^\tau[G,t]$
 over $F^\tau[G']$),
then one can check that  any $x\in F^\tau[G,t]$
satisfies the identity
$$
x^3-3\tr(x)x^2+\bigg(\frac{9}{2}\tr(x)^2-\frac{3}{2}\tr(x^2)\bigg)x
-\tr(x^3)+\frac{9}{2}\tr(x^2)\tr(x)-\frac{9}{2}\tr(x)^3=0.
$$
This guarantees that the bracket is a Lie bracket by
the recognition theorem [BZ, Thm 3.4.7].
We call the Lie algebra $L(\text G_2^{(3)},F^\tau[G,t])$
a {\it twisted $(G, \tau)$-loop algebra} of type $\text G_2^{(3)}$.
We note that if $F$ contains a primitive cubic root of unity,
a twisted $(G, \tau)$-loop algebra of type $\text G_2^{(3)}$ 
can be constructed similarly to the case
of type $\text F_4^{(2)}$.
But our
 $L(\text G_2^{(3)},F^\tau[G,t])$ exists over any $F$.

We often omit the term `untwisted' or `twisted' and
simply say a (locally) $(G, \tau)$-loop algebra.
When $G\cong\Bbb Z$,
we have
$F^\tau[\Bbb Z]\cong F[\Bbb Z]=F[t^{\pm1}]$.
So it is natural to call
 the (locally) $(\Bbb Z, \tau)$-loop algebras above just
{\it (locally) loop algebras},
and of course the loop algebras are
the well-known algebras in Kac-Moody theory.
Also, if $\tau\equiv 1$,
i.e., $F^\tau[G]=F[G]$ is a group algebra,
then
a (locally) $(G, 1)$-loop algebra is simply called
a {\it (locally) $G$-loop algebra}.

If $G$ is a subgroup of $\Bbb Q$,
then $G$ is a directed union of cyclic groups of infinite order, and so
any locally $(G, \tau)$-loop algebra
is a directed union of loop algebras.
Also, if 
$F$ is algebraically closed,
then $F^\tau[G]=F[G]$ by a suitable base change. 
($G$ can be any abelian group 
for this statement, see [P, Lem.2.9] in detail.)

Now,
let $G$ be a subgroup of $\Bbb Q$.
For any two elements $x\otimes t^g$ and  $y\otimes t^h$
in any locally $(G, \tau)$-loop algebra $\Cal L$,
define the new bracket on a 1-dimensional central extension $\tilde{\Cal L}:=\Cal L\oplus Fc$ by
$$[x\otimes t^g,y\otimes t^h]:=[x,y]\otimes \tau(g,h) t^{g+h}+(x,y)\tau(g,h)\delta_{g+h,0}gc
$$
(note $g\in G\subset \Bbb Q\subset F$),
where $(x,y)$ is the trace form or the Killing form depending on the type of $\Cal L$,
or for type $\text B_{\frak I}^{(2)}$,
the direct sum of the trace form and the bilinear form on $V$ determined by 
the symmetric matrix $s$ above,
or for type $\text G_2^{(3)}$,
the direct sum of the trace form on $D_{\Bbb O,\Bbb O}$ and the trace form $t$ 
on $\Bbb O_0$ above.
Indeed, this gives a central extension since
$\tilde{\Cal L}$ is locally an affine Lie algebra,
i.e., a 1-dimensional central extension of a loop algebra,
and $\Cal L$ is a directed union of loop algebras.
One can also show that $\tilde{\Cal L}$ is
a universal central extension of $\Cal L$
[MY2].
Let
$$\hat\Cal L =\tilde{\Cal L}\oplus Fd,$$
where
$d$ is the degree derivation, i.e.,
$$[d,x\otimes t^g]=gx\otimes t^g
\quad\text{and}\quad
[d,c]=0.$$
Let
$$\Cal H =\frak h\oplus Fc\oplus Fd,$$
where
$\frak h$
is the subalgebra of $\Cal L$ consisting of diagonal matrices of degree 0
when $\frak I$ is infinite
or the Cartan subalgebra of 
each finite-dimensional split
simple Lie algebra  $\frak g$
when $\frak I$ is finite.
Then
$\Cal H$ is a Cartan subalgebra of $\hat\Cal L$,
and one can check that
the set of anisotropic roots 
relative to $(\hat\Cal L,\Cal H)$ is a LEARS of null dimension 1.
We also note that $\hat\Cal L$ is an example of 
locally extended affine Lie algebra of null dimension 1 in the sense of [MY1].
In particular, if $G=\Bbb Z$ and $\frak I$ is infinite, then 
the root system of each $\hat\Cal L$
is one of seven reduced LARS listed in Corollary 15,
which is very close to an affine Kac-Moody Lie algebra, 
and we call it
a {\it locally affine Lie algebra}.
In [MY2] we classify locally affine Lie algebras.

\end